\documentclass[sigplan,screen,dvipsnames,nonacm]{acmart} 
\pdfoutput=1

\usepackage{quiver}
\usetikzlibrary{nfold}

\usepackage{amsmath}
\usepackage{amsxtra}
\usepackage{amsthm}

\usepackage{mathrsfs}
\usepackage{tikz-cd}
\usepackage{mathtools}
\usepackage{mathpartir}

\usepackage{graphicx}

\usepackage{todonotes}

\usepackage[super]{nth}

\usepackage[bbgreekl]{mathbbol}
\usepackage{stmaryrd}
\usepackage{enumitem} 
\usepackage{xcolor}
\definecolor{darkgreen}{rgb}{0,0.45,0}
\definecolor{darkred}{rgb}{0.75,0,0}
\definecolor{darkblue}{rgb}{0,0,0.6}
\definecolor{cooltt}{RGB}{19, 92, 183}
\definecolor{redtt}{RGB}{232, 70, 68}
\usepackage{cleveref}

\usepackage{hyperref}
\usepackage[hyphenbreaks]{breakurl}

\usepackage[frozencache]{minted} 

\usetikzlibrary{decorations.markings}
\usetikzlibrary{decorations.pathmorphing}
\tikzset{
  module/.style={
    postaction={decorate},
    decoration={
      markings,
      mark=at position #1 with {\arrow{|}}}},
  module/.default=0.5,
  we/.style=
  { postaction={%
      decorate,
      decoration={
        markings,
        mark=at position #1 with {%
          \node[transform shape, yshift=.2em]{%
            \resizebox{0.5em}{!}{$\sim$}};}}}},
  we/.default=0.5,
  we'/.style=
  { postaction={%
      decorate,
      decoration={
        markings,
        mark=at position #1 with {%
          \node[transform shape, yshift=-.2em, rotate=180]{%
            \resizebox{0.5em}{!}{$\sim$}};}}}},
  we'/.default=0.5,
  iso/.style=
  { postaction={%
      decorate,
      decoration={
        markings,
        mark=at position #1 with {%
          \node[transform shape, yshift=.2em]{%
            \resizebox{0.5em}{!}{$\simeq$}};}}}},
  iso/.default=0.5,
  iso'/.style=
  { postaction={
      decorate,
      decoration={
        markings,
        mark=at position #1 with {%
          \node[transform shape, yshift=-.2em, rotate=180]{%
            \resizebox{0.5em}{!}{$\simeq$}};}}}},
  iso'/.default=0.5,
}

\tikzset{
  proarrow/.style={->, module},
  proequal/.style={-, double, module},
  prodotted/.style={->,dotted, module},
  prodashed/.style={->,dashed, module},
  wearrow/.style={->, we},
  wedashed/.style={->, dashed, we},
  wedotted/.style={->, dotted, we},
  tfibarrow/.style={->>, we=0.45},
  tfibdotted/.style={->>,dotted, we=0.45},
  tfibdashed/.style={->>,dashed, we=0.45},
  tcofarrow/.style={>->, we},
  tcofdashed/.style={>->, dashed, we},
  uwearrow/.style={->, we'},
  uwedashed/.style={->, dashed, we'},
  uwedotted/.style={->, dotted, we'},
  utfibarrow/.style={->>, we'=0.45},
  utfibdotted/.style={->>,dotted, we'=0.45},
  utfibdashed/.style={->>,dashed, we'=0.45},
  utcofarrow/.style={>->, we'},
  isoarrow/.style={->, iso},
  isodashed/.style={->, dashed, iso},
  uisoarrow/.style={->, iso'},
  uisodashed/.style={->, dashed, iso'},
  isocell/.style={=>, iso},
  isocelldashed/.style={=>, dashed, iso},
  uisocell/.style={=>, iso'},
  uisocelldashed/.style={=>, dashed, iso'}
}

\renewcommand{\theenumi}{\roman{enumi}}

\setlist{}
\setenumerate{leftmargin=*,labelindent=0\parindent}
\setitemize{leftmargin=\parindent}

\newtheorem{thm}{Theorem}[section]

\newtheorem{ax}[thm]{Axiom}

\theoremstyle{definition}
\newtheorem{defn}[thm]{Definition}

\theoremstyle{remark}

\makeatletter
\let\c@equation\c@thm
\makeatother
\numberwithin{equation}{section}

\newcommand{\id}{\mathsf{id}}
\newcommand{\ob}{\textup{ob}}

\newcommand{\dhom}{\textup{dhom}}

\newcommand{\cat}[1]{\textup{\textsf{#1}}}

\renewcommand{\AA}{\mathbb{A}}

\DeclareMathAlphabet{\mathbbn}{U}{BOONDOX-ds}{m}{n} 
\newcommand*{\II}{\mathbbn{2}}

\DeclareMathAlphabet{\mathbbe}{U}{bbold}{m}{n}

\newcommand{\Set}{\mathcal{S}\cat{et}}

\makeatletter
\def\makeslashed#1#2#3#4#5{#1{\mathpalette{\sla@{#2}{#3}{#4}}{#5}}}

\def\@mathlower#1#2#3{\setbox0=\hbox{$\m@th#2#3$}\lower#1\ht0\box0}
\def\mathlower#1#2{\mathpalette{\@mathlower{#1}}{#2}}
\makeatother

\newcommand{\univ}{\mathcal{U}}

\newcommand{\refl}{\textup{refl}}

\newcommand{\totalty}[1]{\widetilde{#1}}

\newcommand{\isEquiv}{\mathsf{isEquiv}}
\newcommand{\isContr}{\mathsf{isContr}}
\newcommand{\isProp}{\mathsf{isProp}}
\newcommand{\fib}{\mathsf{fib}}

\newcommand{\evid}{\mathsf{evid}}
\newcommand{\devid}{\mathsf{devid}}
\newcommand{\evrefl}{\mathsf{evrefl}}
\newcommand{\res}{\mathsf{res}}
\newcommand{\unit}{\mathbf{1}}

\newcommand{\tope}{\mathsf{tope}}
\newcommand{\comp}{\mathsf{comp}}
\newcommand{\tyrec}{\mathsf{rec}}
\newcommand{\tyiso}{\mathrm{iso}}
\newcommand{\isIso}{\mathsf{isIso}}
\newcommand{\isSegal}{\mathsf{is}\text{-}\mathsf{pre}\text{-}\infty\text{-}\mathsf{category}}

\newcommand{\isDiscrete}{\mathsf{is}\text{-}\infty\text{-}\mathsf{groupoid}}
\newcommand{\idtoiso}{\mathsf{isoeq}}
\newcommand{\exthtpyeq}{\mathsf{exthtpyeq}}

\newcommand{\idtoarr}{\mathsf{arreq}}

\newcommand{\funext}{\mathsf{funext}}
\newcommand{\extext}{\mathsf{extext}}
\newcommand{\trans}{\mathsf{trans}}
\newcommand{\lift}{\mathsf{lift}}
\newcommand{\htpyeq}{\mathsf{htpy}\text{-}\mathsf{eq}}
\newcommand{\yon}{\mathsf{yon}}
\newcommand{\indpath}{\mathsf{ind}_=}

\newcommand{\exten}[4]{\left\langle\mathchoice{\textstyle\prod_{#1}}{\textstyle\prod_{#1}}{\scriptstyle\prod_{#1}}{\scriptscriptstyle\prod_{#1}} #2 \middle|^{#3}_{#4}\right\rangle}
\newcommand{\ndexten}[4]{\left\langle #1 \to #2 \middle|^{#3}_{#4}\right\rangle}

\def\homtwoshort#1(#2,#3,#4,#5,#6,#7){\hom_{#1}^2(#5,#6;#7)}

\def\homtwoarg#1(#2,#3,#4,#5,#6,#7){\left(
	\begin{tikzpicture}[baseline,xscale=#1]
		\node (a) at (0,0) {$\scriptstyle #2$};
		\node (b) at (1,.3) {$\scriptstyle #3$};
		\node (c) at (2,0) {$\scriptstyle #4$};
		\draw (a) -- node [auto] {$\scriptstyle #5$} (b);
		\draw (b) -- node [auto] {$\scriptstyle #6$} (c);
		\draw (a) -- node [auto,swap] {$\scriptstyle #7$} (c);
	\end{tikzpicture}\right)}

\def\nat#1#2{\underset{#1\to#2}{\hom}}

\newcommand{\rzk}{\textsc{Rzk}}
\newcommand{\Rzk}{\rzk{}}

\newcommand{\libUniMath}{\href{https://github.com/UniMath/UniMath}{\texttt{UniMath}}}
\newcommand{\libagdaunimath}{\href{https://unimath.github.io/agda-unimath/}{\texttt{agda-unimath}}}
\newcommand{\libagdacategories}{\href{https://github.com/agda/agda-categories}{\texttt{agda-categories}}}
\newcommand{\libmathlib}{\href{https://github.com/leanprover-community/mathlib}{\texttt{mathlib}}}

\makeatletter
\newcommand{\yonedawebsiteref}[2]{\if@ACM@anonymous\href{https://anonymous-yoneda.github.io/yoneda/#1}{#2}%
\else\href{https://emilyriehl.github.io/yoneda/CPP-2024/#1}{\small\color{OliveGreen}{#2}}%
\fi
}
\makeatother
\newcommand{\inlinerzk}[1]{\mintinline{rzk}{#1}}
\newcommand{\rzkref}[2]{\yonedawebsiteref{#1\#define:#2}{\texttt{#2}}}
\newcommand{\rzkrefU}[3]{\yonedawebsiteref{#1\#define:#2}{\texttt{#3}}}
\newcommand{\rzkhottsecref}[1]{\yonedawebsiteref{hott/#1.rzk}{\texttt{#1}}} 
\newcommand{\rzksimpsecref}[1]{\yonedawebsiteref{simplicial-hott/#1.rzk}{\texttt{#1}}} 

\DeclareUnicodeCharacter{0394}{\ensuremath{\Delta}}   
\DeclareUnicodeCharacter{2261}{\ensuremath{\equiv}}   
\DeclareUnicodeCharacter{2081}{\ensuremath{_1}}       
\DeclareUnicodeCharacter{2082}{\ensuremath{_2}}       
\DeclareUnicodeCharacter{21A6}{\ensuremath{\mapsto}}  
\DeclareUnicodeCharacter{2264}{\ensuremath{\leq}}     
\DeclareUnicodeCharacter{03A3}{\ensuremath{\Sigma}}   
\DeclareUnicodeCharacter{03C6}{\ensuremath{\phi}}     
\DeclareUnicodeCharacter{03C8}{\ensuremath{\psi}}     
\DeclareUnicodeCharacter{2227}{\ensuremath{\land}}    
\DeclareUnicodeCharacter{22A5}{\ensuremath{\bot}}     
\DeclareUnicodeCharacter{221E}{\ensuremath{\infty}}   

\newcommand{\pair}[2]{(#1,#2)}

\let\jdeq\equiv
\newcommand*{\defeq}{\mathrel{\vcentcolon\jdeq}}

\newcommand{\ie}{\emph{i.e.}}
\newcommand{\aka}{\emph{aka}}
\newcommand{\eg}{e.g.}
\newcommand{\cf}{cf.}

\setcopyright{rightsretained}
\acmDOI{10.1145/3636501.3636945}
\acmYear{2024}
\copyrightyear{2024}
\acmSubmissionID{poplws24cppmain-p18-p}
\acmISBN{979-8-4007-0488-8/24/01}
\acmConference[CPP '24]{Proceedings of the 13th ACM SIGPLAN International Conference on Certified Programs and Proofs}{January 15--16, 2024}{London, UK}
\acmBooktitle{Proceedings of the 13th ACM SIGPLAN International Conference on Certified Programs and Proofs (CPP '24), January 15--16, 2024, London, UK}
\received{2023-09-19}
\received[accepted]{2023-11-25}

\begin{document}

\title{Formalizing the \texorpdfstring{$\infty$}{∞}-Categorical Yoneda Lemma}

\author{Nikolai Kudasov}
\email{n.kudasov@innopolis.ru}
\orcid{0000-0001-6572-7292}
\affiliation{%
  \institution{Lab of Programming Languages and Compilers\\Innopolis University}
  \streetaddress{Universitetskaya 1}
  \city{Innopolis}
  \state{Tatarstan Republic}
  \country{Russia}
  \postcode{420500}
}

\author{Emily Riehl}
\orcid{0000-0002-8465-8859}
\authornote{The second and third authors are supported by the US Army Research Office under MURI Grant W911NF-20-1-0082. In addition, second author is supported by the US National Science Foundation via the grant DMS-2204304, by the US Air Force Office of Scientific Research under award number FA9550-21-1-0009, and by a Simons Fellowship with award number 920415.}
\email{eriehl@jhu.edu}
\author{Jonathan Weinberger}
\orcid{0000-0003-4701-3207}
\authornotemark[1]
\email{jweinb20@jhu.edu}
\affiliation{%
  \institution{Department of Mathematics\\Johns Hopkins University}
  \streetaddress{3400 N Charles Street}
  \city{Baltimore}
  \state{MD}
  \country{USA}
  \postcode{21218}}


\begin{abstract}
  Formalized $1$-category theory forms a core component of various libraries of mathematical proofs. However, more sophisticated results in fields from algebraic topology to theoretical physics, where objects have ``higher structure,'' rely on infinite-dim\-en\-sional categories in place of $1$-dimensional categories, and $\infty$-cat\-e\-gory theory has thusfar proved unamenable to computer formalization.

	Using a new proof assistant called \rzk,
which is designed to support Riehl--Shulman's  simplicial extension of homotopy type theory for synthetic $\infty$-category theory, we provide the first formalizations of results from $\infty$-category theory. This includes in particular a formalization of the Yoneda lemma, often regarded as the fundamental theorem of category theory, a theorem which roughly states that an object of a given category is determined by its relationship to all of the other objects of the category. A key feature of our framework is that, thanks to the synthetic theory, many constructions are automatically natural or functorial.
  We plan to use \rzk{} to formalize further results from $\infty$-category theory, such as the theory of limits and colimits and adjunctions.

	  \end{abstract}



\maketitle


\section{Introduction}\label{sec:intro}

\emph{Computer proof assistants} are computer programs that formally verify the logical reasoning of mathematical proofs. There are a wide variety of such programs --- including Agda \cite{BoveDybjerNorell2009,NorellChapman2009}, Coq \cite{BertotCasteran2013}, HOL~Light \cite{Harrison2009}, Isabelle \cite{NipkoWenzelPaulson2002}, and Lean \cite{deMouraKongAvigadVanDoornVonRaumer2015,deMouraUllrich2021}, among many others --- which have enjoyed striking success in recent decades. Notable accomplishments include:

\begin{itemize}
\item  A project \cite{HalesAdamsBauerDangHarrisonLeTurongKaliszykMagronMacLaughlinNguyenEtAl2017} in HOL~Light from 2003--2014, to formally verify the Ferguson--Hales proof of the Kepler conjecture, after the referees for the 2005 \emph{Annals} publication issued a disclaimer stating that they were only 99\% certain of its correctness.
\item A project \cite{GonthierAspertiAvigadBertotCohenGarillotLeRouxMahboubiOConnorOuldBiha2013} in Coq, from 2006--2012, to formally verify the Feit--Thompson Odd Order Theorem, a foundational result in the classification of finite simple groups.
\item A project in Lean \cite{deMouraUllrich2021}, from 2020--2022, dubbed the ``liquid tensor experiment'' \cite{Scholze2022}, to formally verify a result from condensed mathematics after Peter Scholze expressed concern about the correctness of his own proof.
\end{itemize}

Part of the task in formalizing a cutting-edge mathematical result is to develop an accompanying \emph{library} of background mathematics on which it depends. For instance, the liquid tensor experiment required a formalized library of standard results from graduate-level homological algebra as well as many other background topics.

In addition to homological algebra, Lean's mathematics library \libmathlib{}~\cite{mathlib2020} contains standard results from number theory, representation theory, general topology, linear algebra including Banach and Hilbert spaces, measures and integral calculus, random variables, basic algebraic geometry, model theory, and category theory, among other topics. Despite all these achievements, as of the writing of this article, Lean's \libmathlib{} does not contain any $\infty$-category theory, and thus myriad recent results from algebraic K-theory \cite{BlumbergGepnerTabuada2013}, derived and spectral algebraic geometry \cite{Lurie2004,Lurie2018}, the Langlands program \cite{FarguesScholze2021}, and symplectic geometry \cite{NadlerTanaka2020} are inaccessible to formalization.

Building such a library would arguably be no more difficult than the successes mentioned above, but the endeavor would be frustrated by the ways in which the traditional set-based foundations for mathematics are not optimal for reasoning about $\infty$-cat\-e\-go\-ries. To give a precise definition of an $\infty$-\emph{category} --- which, roughly speaking, is an infinite-dimensional category with a weak composition law in which all morphisms above dimension 1 are weakly invertible --- one must pick a ``model,'' a Bourbaki-style mathematical structure presenting the $\infty$-categorical data \cite{Camarena2013,Bergner2018}. A variety of models --- such as \emph{quasi-categories} \cite{BoardmanVogt1973,Joyal2002}, \emph{complete Segal spaces} \cite{Rezk2001}, and \emph{Segal categories} \cite{HirschowitzSimpson2001,Pellissier2002} --- are used at various places in the literature, and theorems are often proven ``analytically,'' in reference to the ``coordinates'' of a particular model. A computer formalizer is thus faced with an unattractive choice of either
\begin{itemize}
    \item  picking one model, which must then be used for the entire library of subsequent results, or
    \item formalizing multiple models and the comparisons between them \cite{JoyalTierney2007}, which significantly increases the workload.\footnote{Experts in the field often prefer to work ``model-independently'' \cite{AyalaMazelGeeRozenblyum2022,Lurie2003} which can be done either by using $\infty$-category theory itself as the ambient metatheory, or deploying the formalism of \emph{$\infty$-cosmoi} (\ie, categories of $\infty$-categories) \cite{RiehlVerity2022}, but either approach would require some initial formalization in a specific model of $\infty$-categories.}
\end{itemize}

\subsection{Reimagining the foundations of \texorpdfstring{$\infty$}{infinity}-category theory}

A radical-sounding alternative, which we argue is worth taking seriously, is to change the foundation system. The article ``Could $\infty$-category theory be taught to undergraduates?'' \cite{Riehl2023}  argues that it is possible to narrow the gap between $\infty$-category theory and ordinary 1-category theory by replacing the traditional foundations with a directed extension of \emph{homotopy type theory} \cite{Rijke2022,hottbook}. The basis for this claim is the paper \cite{RiehlShulman2017} (and the follow-up work of \cite{BardomianoMartinez2022,BuchholtzWeinberger2023,Weinberger2022PhD}), which develops the basic theory of $\infty$-categories in an alternative foundational framework established there. The \emph{simplicial type theory} is a formal framework that permits one to make the following intuitive definitions rigorous:
\begin{itemize}
    \item A type is a \emph{pre-$\infty$-category} (aka a \emph{Segal type}) if every composable pair of arrows has a unique composite.
    \item A pre-$\infty$-category is an \emph{$\infty$-category} (aka a \emph{Rezk type}) if equalities are equivalent to isomorphisms.
    \item A type is an \emph{$\infty$-groupoid} (aka a \emph{discrete type}) if equalities are equivalent to arrows.
    \item A type family is a \emph{covariant fibration} (aka a \emph{covariant type family}) if every arrow in the base type has a unique lift with specified domain.
\end{itemize}

The intended model of this formal system is in the category of simplicial spaces (aka the category of bisimplicial sets), which Shulman had shown provides a model of homotopy theory, in which types are interpreted as Reedy fibrant simplicial spaces \cite{Shulman2015}. In this model, pre-$\infty$-categories correspond to Segal spaces \cite{Segal1968,Rezk2001}, $\infty$-categories correspond to complete Segal spaces \cite{Rezk2001}, and covariant fibrations correspond to left fibrations \cite{deBrito2016,KazhdanVarshavsky2014}. The phrases ``for all \ldots there exists \ldots unique'' are meant in the standard sense of homotopy type theory \cite{Rijke2022,hottbook}. In particular, following the homotopical extension of the Curry-Howard correspondence \cite{HofmannStreicher1998,AwodeyWarren2009,Voevodsky2009}, \emph{uniqueness} means \emph{contractibility} --- which is precisely what is true semantically for the composition operation in an $\infty$-category.\footnote{Those familiar with the Segal space model of $\infty$-categories may be surprised that the definition of a pre-$\infty$-category refers to binary sequences of composable arrows and not also composable triples and quadruples and so on. Here the binary statement subsumes the $n$-ary ones for $n \geq 0$ because it is interpreted \emph{internally} in the model as the assertion that the internal mapping types mapping out of the $2$-simplex and out of its inner horn are equivalent~\cite[Section~5]{RiehlShulman2017}.} This model validates the terminology used in our synthetic theory of $\infty$-categories. Via this interpretation functor, proofs of theorems about $\infty$-categories in the simplicial type theory really do prove theorems about $\infty$-categories, as instantiated by the complete Segal spaces model in traditional foundations.

More generally, Shulman has proven that homotopy type theory has semantics in any $\infty$-topos \cite{Shulman2019,Riehl2022} and Weinberger \cite{Weinberger2022} has shown that the simplicial type theory of \cite{RiehlShulman2017} can be interpreted in simplicial objects in any $\infty$-topos. Thus, theorems proven about synthetic $\infty$-categories in the simplicial type theory also apply to the \emph{internal $\infty$-categories} studied by Martini, Rasekh, Stenzel, and Wolf among others \cite{Rasekh2022,Stenzel2023,MartiniWolf2021,Martini2022}.

\subsection{Formalizing \texorpdfstring{$\infty$}{∞}-category theory}

It is relatively standard practice in homotopy type theory to formalize results while writing the corresponding paper proofs.\footnote{While homotopy type theory cannot be formalized in Lean or Idris~\cite{Brady2013} because their kernels assume that all types are sets, contradicting Voevodsky's \emph{univalence axiom},
it can be done in Agda, Coq, and a growing variety of experimental proof assistants.}
At the time of the writing of the first paper on synthetic $\infty$-category theory by Riehl and Shulman \cite{RiehlShulman2017},
it was not possible to formalize any of its results because the work is done in an \emph{extension} of traditional homotopy type theory,
with multilevel contexts and a new type-forming operation providing \emph{extension types}. 

The multilevel context includes \emph{cube} and \emph{tope} layers, upon which the final layer of types may depend. The cube and tope layers combine to provide directed \emph{shapes} that parametrize arrows, composable pairs of arrows, and their composites, among other \emph{simplices} and their \emph{subshapes}. The rules for these layers are spelled out in \cite[\S 2.1]{RiehlShulman2017}. The \emph{extension types}
reify all the possible totalizations or \emph{extensions} of a given partially defined datum along a subshape inclusion and satisfy rules enumerated in \cite[\S 2.2]{RiehlShulman2017}.

Extension types are playing an important rule in both simplicial~\cite{RiehlShulman2017} and cubical type theory~\cite{BezemCoquandHuber2014,CCHM2018,OrtonPitts2018,Awodey2018,AngiuliBrunerieCoquandHarperHouLicata2021,CMS2020}, the metatheory thereof~\cite{Sterling2021,SterlingAngiuli2021,Sterling2023,Gratzer2022,AagaardKristensenGratzerBirkedal2022,WABN2022}, and various applications in programming language theory~\cite{GSACB2022,SterlingHarper2021,Zhang2023}. In this work, we require at least \emph{simplicial extension types} --- where the adjective ``simplicial'' refers to the geometry of the subshape inclusions --- whereas cubical provers, such as Cubical Agda~\cite{VezzosiMoertbergAbel2021}, red* family of proof assistants~\cite{RedPRL,redtt,cooltt}, Aya, and Arend all support only \emph{cubical extension types}. A reasonable alternative could be to postulate extension types, e.g. in Agda with user-defined rewrite rules.\footnote{see \url{https://agda.readthedocs.io/en/latest/language/rewriting.html}}
However, to the best of our knowledge, this approach fails to capture rewriting for an application of a variable of extension type to an argument,
leaving computational rules incomplete, and requiring additional bookkeeping in user-defined proofs to push computation further.
The new proof assistant \rzk{} developed by Kudasov~\cite{Kudasov2023rzk} supports \emph{simplicial extension types},\footnote{Technically, \rzk{} supports extension types over arbitrary topes, except only simplicial and trivial (unit) topes are available in \rzk{} v0.6.7.} and, thus, finally one can formally test the claims made in the article \cite{Riehl2023}. This is the content of our project.

In \S\ref{sec:stt}, we describe the simplicial type theory, and in \S\ref{sec:syn-inftycat} we introduce synthetic $\infty$-category theory. In \S\ref{sec:rzk}, we describe the \Rzk~proof assistant. In \S\ref{sec:yoneda}, we describe our formalization of the $\infty$-categorical Yoneda lemma in \Rzk. In \S\ref{sec:comparison}, we compare this formalization with parallel formalizations of the $1$-categorical Yoneda lemma in both traditional and univalent foundations. In \S\ref{sec:future}, we offer a few takeaways from this formalization project and describe related future work.

This paper is accompanied by an open source formalization repository available at \url{https://github.com/emilyriehl/yoneda}.
The code is written in a literate style, using \rzk{} code blocks embedded in Markdown files.
For the convenience of the readers, formalizations are rendered and available at \url{https://emilyriehl.github.io/yoneda/}.
In this paper, most definitions and proofs have a corresponding identifier, acting also as a hyperlink to its formalized counterpart.
For example, the Yoneda lemma is formalized as \rzkref{simplicial-hott/09-yoneda.rzk}{yoneda-lemma}.

\subsection{Contributions}

Our contribution consists of a library for synthetic $\infty$-category theory, developed from scratch for the new \Rzk\, proof assistant. The library contains large portions of formalized synthetic $\infty$-category theory from the paper \cite{RiehlShulman2017}, previously not supported by any proof assistant. The results range from properties of extension types, to formal properties of synthetic $\infty$-categories and fibrations thereof, to the Yoneda lemma. We also formalized many results from Standard/Book HoTT~\cite{hottbook,Rijke2022} that were needed as a foundation. Finally, we contribute a comparison to other formalizations of the Yoneda lemma in other systems and proof assistants, in particular to the formalization of the Yoneda lemma for precategories that we contributed to the \libagdaunimath{} library.

Moreover, the formalization process led us to discover a mistake in the paper \cite{RiehlShulman2017}: the published proof of the ``only if'' direction of Proposition 8.13 employed circular reasoning.\footnote{While this mistake could have been caught by the original authors, the process of formalization made it entirely transparent: a proof of the conclusion was called as a hypothesis by one of the lemmas cited in the original argument.} Fortunately, the stated result remains true. Our new formalized proof (\rzkref{simplicial-hott/08-covariant.rzk}{is-segal-is-covariant-representable}) now appears in \cite{RiehlShulman2023arXiv-v5}. 

In parallel with the submission of this paper, we  invited other researchers to contribute to the broader project of formalizing synthetic $\infty$-category theory, extending the results reported upon here. To that end, we created a clone of our repository,\footnote{\url{https://github.com/rzk-lang/sHoTT}} which to date has had a dozen contributors, in addition to ourselves. In that repository, together with our new collaborators, we are already pursuing some of the projects we describe in \S\ref{sec:future}.

\subsection{Prerequisites}

Our work involves a subtle interplay between  (homotopy) type theory and (higher) category theory. While \S\ref{sec:stt} and \S\ref{sec:syn-inftycat} provide some introductory discussion of the relevant concepts, some familiarity with these topics is necessary for a deeper understanding of our work, and thus we  direct the interested reader to further background sources in the literature. Our work assumes very solid knowledge of Martin-Löf type theory (MLTT) and homotopy type theory (HoTT); see the two books on homotopy type theory due to the Univalent Foundations Project~\cite{hottbook} and Rijke~\cite{Rijke2022}, resp. A working mathematician's discussion of $\infty$-categories and their peculiarities  is given by Riehl in~\cite{Riehl2023}. For the discussion of the \Rzk~proof assistant and the comparisons of formalizations, it may be helpful to have some familiarity with Coq~\cite{BertotCasteran2013} and (Cubical) Agda~\cite{BoveDybjerNorell2009,NorellChapman2009,VezzosiMoertbergAbel2021}, as well as \libagdaunimath{}~\cite{UniMath}.

\subsection{Related work}

A roughly parallel synthetic framework for $\infty$-category theory has been proposed by Weaver and Licata using bicubical sets as an intended model \cite{WeaverLicata2020}. An alternate approach to formalizing higher category is within the framework of \emph{two-level type theory}, using extensional type theory as a meta-theory, see \eg~\cite{Voevodsky2013,AnnenkovCapriottiKrausSattler23,Kraus2021}.

A conceptual discussion of the approach behind simplicial type theory with comparisons is done by Buchholtz in~\cite{Buchholtz2019}. A self-contained overview of both syntactic and semantic aspects of simplicial type theory is given in the master's thesis of Bakke~\cite{Bakke2021}.

Furthermore, there has been extensive work on directed type theories \cite{LicataHarper2011,Warren2013,Nuyts2015,North2019}, though most of this was not created to describe $\infty$-category theory. Other work includes domain-specific languages for two-dimensional categories of various flavors, see~\cite{Garner2009,AhrensNorthvdWeide2022}, as well as further work for the case of bicategories~\cite{Mohri1997,AhrensMaggesi2018,FioreSaville2019,Stark2023,AhrensNorthvdWeide2023}, virtual equipments~\cite{NewLicata2023}, and double categories~\cite{vdWRAN2023}. There also exist other type theories capturing infinite-dimensional categorical structures. Notable developments include~\cite{FinsterRiceVicary2021,FinsterReutterVicaryRice2022,Dean2020,FinsterMimram2017,BenjaminFinsterMimram2021,AlliouxFinsterSozeau2021,BordgDonaMateo2023}.
However, these systems differ from the one that we are using in two major aspects: their setup and their purposes. Our framework features a synthetic and homotopical theory of $\infty$-categories with the aim of developing a range of classical $\infty$-categorical results. The other frameworks tend to involve a specific model of either strict or weak infinite-dimensional categories.

Aside from direct applications to category theory, new kinds of type theories have been devised for the purpose of doing differential topology and stable homotopy theory  synthetically, making heavy use of type-theoretic \emph{modalities}~\cite{Shulman2018,SchreiberShulman2014,Cherubini2022,RijkeShulmanSpitters2020,Myers2022Fib,Myers2021,Myers2022Orbi,RileyFinsterLicata2021,MyersRiley2023}.

\section{The Simplicial Type Theory}\label{sec:stt}

In~\cite{RiehlShulman2017}, Riehl--Shulman develop a type theory to reason synthetically about $\infty$-categories. The key features of their theory is that $\infty$-categories can be described in relatively simple terms, and all the results are invariant under \emph{homotopy equivalence} --- the right notion of equivalence of $\infty$-categories. This is in stark contrast to the more traditional and familiar developments of $\infty$-category theory in set theory, \cf~\eg~\cite{Lurie2009,Joyal2008}. We will give an overview of the structure and features of the simplicial type theory, with an emphasis on its use for synthetic $\infty$-category theory.

Simplicial type theory shares some concepts with cubical type theory~\cite{BezemCoquandHuber2014,CCHM2018,OrtonPitts2018,Awodey2018,AngiuliBrunerieCoquandHarperHouLicata2021,CMS2020}. A comparison, at least of the underlying shape theories, is provided by Riehl--Shulman in~\cite[3.2]{RiehlShulman2017}. Indeed, Weaver--Licata develop concepts from the paper \cite{RiehlShulman2017} in a (bi)cubical directed type theory~\cite{WeaverLicata2020}.

The theory builds on Martin-L\"{o}f intensional type theory (MLTT)~\cite{MartinLoef1975} whose \emph{intensional identity types} have homotopically well-behaved \emph{path objects} as models~\cite{AwodeyWarren2009,KapulkinLumsdaine2021,Streicher2014,Riehl2023-HoTT}. This homotopical interpretation, paired with Voevodksy's \emph{univalence axiom}, which allows one to treat homotopy equivalent types as (intensionally) equal, goes by the name \emph{homotopy type theory} (HoTT) or \emph{univalent foundations}   \cf~\cite{Voevodsky2009,AwodeyWarren2009,hottbook}. Homotopy type theory may be thought of as a synthetic theory for $\infty$-groupoids (aka homotopy types) and thus provides a fertile basis for the simplicial type theory.

\subsection{Base Theory: Martin-L\"{o}f Intensional Type Theory}\label{ssec:mltt}

\subsubsection*{Overview}

The base theory is intensional Martin-L\"{o}f type theory~\cite{MartinLoef1975} with $\Sigma$-, $\Pi$-, and identity types. Though $\rzk$ works with a universe type to implement dependent types, this assumption is not necessary~\cite[Remark~2.5]{RiehlShulman2017}.\footnote{In particular, though convenient for certain applications, univalence is not necessary for our development.} To stay in line with the notation of~\cite{RiehlShulman2017}, we also notate a dependent type $x:A \vdash C(x)$ as a \emph{type family} $C:A \to \univ$, pretending $\univ$ is a universe type (without being explicit about universe hierarchy or different levels of size).

\subsubsection*{$\Sigma$-Types (\rzkhottsecref{05-sigma})}

The type formers $\Sigma$ and $\Pi$, resp., generalize existential and universal quantification, resp., as follows. For $C: A \to \univ$, the \emph{dependent sum} $\sum_{x:A} C(x)$ is the type consisting of dependent pairs $\pair{a}{c}$ with $a:A$ and $c:C(a)$. This is also referred to as the \emph{total type} of the family $C$. The $\Sigma$-type comes with the usual set of rules for formation, introduction (by forming dependent pairs), and elimination (by projecting to the factors). We also assume the $\beta$- and $\eta$-computation rules to be satisfied, meaning that introduction and elimination are inverse to each other in the strictest possible way, \ie, \emph{up to judgmental equality}.

The family $C:A \to \univ$ can alternatively be encoded as a map $p_C : \totalty{C} \to A$, with the total type $\totalty{C} :\jdeq \sum_{x:A} C(x)$, and the projection $p_C(a,c) :\jdeq a$. The total type is then the ``sum'' of all the fibers of $C$, canonically indexed by $A$. If $C$ is a constant family, \ie, $C(a) \jdeq B$ for all $a:A$ and some type $B$, the $\Sigma$-type becomes the \emph{cartesian product} $A \times B$.

\subsubsection*{$\Pi$-Types}

Of particular interest is the notion of \emph{dependent function} or \emph{section} of a family $C: A \to \univ$, which is an assignment $\sigma$ to each element $x:A$ of some element $\sigma(x) : C(x)$ in the corresponding fiber. This is reified as the \emph{dependent product} type $\prod_{x:A} C(x)$, with introduction rule given by $\lambda$-abstraction and elimination rule by function application. Likewise, we require the $\beta$- and $\eta$-rules to hold judgmentally. When the type family $C$ is constant with value some type $B$, the dependent function type reduces to an ordinary function type, denoted by $A \to B$ or $B^A$.

\subsubsection*{Identity Types (\rzkhottsecref{01-paths})}\label{sssec:id}

The \emph{Martin-L\"{o}f identity types} $(a=_Ab)$ for a type $A$ and elements $a,b:A$ capture the idea that equality between terms of a type is witnessed proof-relevantly by a term $p:a=_Ab$. In the homotopical models, identity types get interpreted as \emph{path objects} in the sense of homotopical algebra~\cite{AwodeyWarren2009}, so elements $p:(a=_Ab)$ can be seen as paths from $a$ to $b$ in $A$. The introduction rule is given by the canonical \emph{reflexivity terms} $\refl_a : (a=_Aa)$ witnessing self-identity. Elimination is given by the \emph{path induction principle}. Intuitively, this says the following. First, for a type $A$, fix $a:A$. Then, for a family $C:\left( \sum_{x:A} (a =_A x)\right) \to \univ$ the type of sections $\prod_{\pair{y}{p} : \sum_{x:A}(a=_Ax)} C(y,p)$ is equivalent to $C(a,\refl_a)$ via the map
\[
	\evrefl_{C,a} : \left( \prod_{\pair{y}{p} : \sum_{x:A}(a=_Ax)} C(y,p) \right) \to C(a,\refl_a).
\]
In particular, given $d:C(a,\refl_a)$ we obtain a section
\[
	\indpath(d) : \prod_{\pair{y}{p} : \sum_{x:A}(a=_Ax)} C(y,p) \tag{\rzkref{hott/01-paths.rzk}{ind-path}}
\]
such that $\indpath(d)(a,\refl_a) \jdeq d$. Thus, for type families over (based) path types, to produce a section of the whole family it suffices to produce a section only at the reflexivity loop.

\subsubsection*{The Homotopy Theory of Types}\label{sssec:hott-types}

The following notions are due to Voevodsky~\cite{Voevodsky2009}, \cf~also \cite{AwodeyWarren2009,KapulkinLumsdaine2021,Streicher2014,Riehl2023-HoTT}.
According to the idea that terms $p:(a=_Ab)$ encode paths in a type we want to express when a type is homotopically trivial \aka~\emph{contractible}. This is witnessed by the type
\[ \isContr(A) \defeq \sum_{x:A} \prod_{y:A} (x=_A y). \tag{\rzkref{hott/06-contractible.rzk}{is-contr}}\]
 A contractible type $A$ comes equipped with a canonical inhabitant, the \emph{center of contraction} $c_A : A$  (\rzkref{hott/06-contractible.rzk}{center-contraction}) and a homotopy $H_A : \prod_{y:A} (c_A =_A y)$ (\rzkref{hott/06-contractible.rzk}{homotopy-contraction}). Contractible types are equivalent to the point or \emph{terminal type} $\unit$, see (\rzkref{hott/06-contractible.rzk}{contr-iff-terminal-map-is-equiv}).

Traditional homotopy theory involves constructions on topological spaces that are invariant under \emph{homotopy equivalence}, which is a pair of maps between two spaces in opposite directions whose composites are homotopic to the identity. Translating this into type theory, a map $f:A \to B$ between types is a (homotopy) equivalence when there is a term inhabiting the type
\[ \isEquiv(f) \defeq \sum_{g:B \to A} (g \circ f =_{A\to A} \id_A) \times  \sum_{h:B \to A} (f \circ h =_{B \to B} \id_B). \tag{\rzkref{hott/03-equivalences.rzk}{is-equiv}}\]
This type is a \emph{proposition} in the sense that it is contractible whenever if it is inhabited. By \cite[12.1.3]{Rijke2022}, this can equivalently be captured by the type
\[ \isProp(A) \defeq \prod_{x,y:A} (x=_Ay). \tag{\rzkref{hott/09-propositions.rzk}{is-prop}}\]
When a type $A$ is a proposition (\ie, $\isProp(A)$ is inhabited), then it can be treated as a \emph{mere property} (up to homotopy, \ie, a contractible and thus trivial choice of data) rather than additional \emph{structure}. The fact that $\isEquiv(f)$ is always a proposition hence means that being an equivalence is, in fact, a \emph{property} of a map, much in line with the expected intuition. It turns out there is a further equivalent characterization of when a map is an equivalence in that sense, namely if and only if all its \emph{fibers} $\fib(f,b) \defeq \sum_{x:A} (f(x) =_A b)$ are contractible, \ie,
\[ \isEquiv(f) \simeq \prod_{b:B} \isContr(\fib(f,b)). \]
If type families are understood as \emph{fibrations} $p: \sum_{x:A} C(x) \to A$, then equivalences in this sense behave like \emph{trivial fibrations} (\rzkhottsecref{10-trivial-fibrations}) whose fibers are all contractible. These homotopical interpretations of Martin L\"{o}f's dependent type theory open up a whole area of research doing homotopy theory synthetically, \cf~\cite{hottbook,Rijke2022}.

\subsubsection*{Function Extensionality (\rzkref{hott/03-equivalences.rzk}{FunExt})}\label{sssec:funext}

While we do not require the univalence axiom in our formalization, we do make use of \emph{function extensionality}, which is one of its consequences \cite[Theorem~17.3.2]{Rijke2022}: we will postulate the map
\[ \htpyeq :  \prod_{X:\univ} \prod_{A: X \to \univ} \prod_{f,g: \prod_X  A} (f = g) \to \prod_{x:X} (fx = gx)\]
defined via path induction by
\[ \htpyeq(X,A,f,f,\refl_f,x)\defeq \refl_{fx} \tag{\rzkref{hott/03-equivalences.rzk}{htpy-eq}} \]
is an equivalence, \ie, there exists a term
\[ \funext : \prod_{X:\univ}  \prod_{A: X \to \univ} \prod_{f,g: \prod_X  A}  \isEquiv(\htpyeq_{X,A,f,g}). \]

\subsubsection*{The $\infty$-Groupoid Structure on a Type}

By (iterated) path induction one can prove the existence of functions
\begin{align*}
	&\prod_{x,y:A} (x=_A y) \to (y=_A x), \tag{\rzkref{hott/01-paths.rzk}{rev}} \\
 	&\prod_{x,y,z:A} (x =_A y) \to (y =_A z) \to (x =_A z) \tag{\rzkref{hott/01-paths.rzk}{concat}}
\end{align*}  serving to reverse paths as well as concatenating them. One can show that these satisfy the expected groupoid laws, but only up to propositional equality, endowing every type canonically with the structure of a (weak) $\infty$-groupoid, \cf~\cite{HofmannStreicher1998,vandenBergGarner2011}.

While $\infty$-groupoids are special cases of $\infty$-categories, in a general $\infty$-category we require directed ``arrows'' that are not necessarily reversible. This suggests the following extensions of the underlying type theory.

\subsection{Extension~1: Cube and Tope Layers}\label{ssec:topes}

Intuitively, a synthetic $\infty$-category is a type where directed arrows can be composed up to homotopy. To reason about directed arrows, their composites, and other shapes arising from this the idea is to introduce an appropriate shape theory to the type theory. The shapes will be part of the contexts so that type families and sections can depend on them.

Each shape is viewed as a \emph{subshape} embedded inside a \emph{higher dimensional (directed) cube}. This is reminiscent of the basic setup of \emph{cubical type theory}.

For the \emph{cube layer}, consider a new pretype $\II$, equipped with two distinct elements $0,1:\II$, and a binary relation $\le$ making $\II$ into a strict partial order with bottom element $0$ and top element $1$. The Lawvere theory generated by $\II$ constitutes the cube layer, \ie, the cubes are exactly the finite powers $\II^n$, with $\II^0 \jdeq \unit$. The partial order is captured by a new judgment form, called a tope:
\[ x,y:\II \vdash x \le y \, \tope \]

The \emph{tope layer} is a finitary intuitionistic logic over the cube layer. The intention is to carve out \emph{subshapes} $\Phi \subseteq I$ of a cube $I$ by describing it via a formula on the cube variables. In general: if $I$ is a cube and $\varphi$ is a tope in context $t:I$, written as a judgment $t:I \vdash \varphi \, \tope$, then $\Phi \defeq \{ t :  I \; | \; \varphi \}$ is the  \emph{shape} corresponding to $\varphi$. This way, one can define important shapes such as the $n$-simplex $\Delta^n$, for $n \in \mathbb N$, its boundaries $\partial \Delta^n$, the \emph{$(n,k)$-horns} $\Lambda_k^n$ for $k \le n$, and more. E.g., we have the following formulas, \cf~also~\Cref{fig:shapes}:
\begin{align*}
	\Delta^1 & \defeq \{t : \II \mid \top\} \subseteq \II \\
	\partial\Delta^1 & \defeq \{t : \II \mid (t \jdeq 0) \lor (t \jdeq 1)\} \subseteq \II \\
	\Delta^2 & \defeq \{ \pair{t}{s} : \II^2 \mid s \leq t \} \subseteq \II^2 \\
	\partial\Delta^2 & \defeq \{ \pair{t}{s} : \II^2 \mid (s \jdeq 0) \lor (s \jdeq t) \lor (t \jdeq 1) \} \subseteq \II^2\\
	\Lambda_1^2 & \defeq \{ \pair{t}{s} : \II^2 \mid (s \jdeq 0) \lor (t \jdeq 1) \} \subseteq \II^2
	\tag{\rzksimpsecref{03-simplicial-type-theory}}
\end{align*}

Like in cubical type theory, we connect the standard type layer with the cube and tope layer through a three-part context, which allows type families $A$ to depend on a cube context $\Xi$, a tope context $\Phi$, and a type context $\Gamma$, written as $\Xi \mid \Phi  \mid \Gamma \vdash A$.

The directed arrows in a type are now defined using our interval shape $\Delta^1$ and another feature to be introduced, the \emph{extension types}.

{\footnotesize
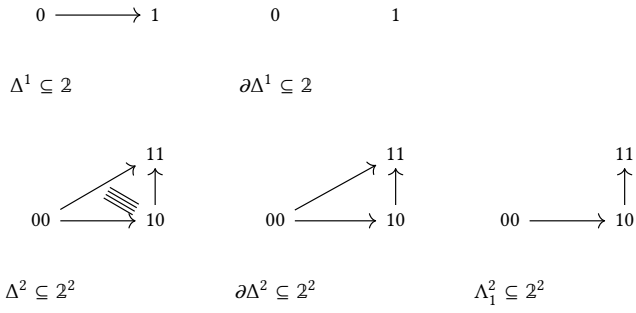
\begin{figure}
	\[\begin{tikzcd}
		0 & 1 & 0 & 1 \\
		{\Delta^1 \subseteq \II} && {\partial \Delta^1 \subseteq \II} \\
		& 11 && 11 && 11 \\
		00 & 10 & 00 & 10 & 00 & 10 \\
		{ \Delta^2 \subseteq \II^2} && {\partial \Delta^2\subseteq \II^2} && {\Lambda_1^2 \subseteq \II^2}
		\arrow[from=1-1, to=1-2]
		\arrow[from=4-2, to=3-2]
		\arrow[""{name=0, anchor=center, inner sep=0}, from=4-1, to=3-2]
		\arrow[from=4-3, to=4-4]
		\arrow[from=4-4, to=3-4]
		\arrow[from=4-1, to=4-2]
		\arrow[from=4-3, to=3-4]
		\arrow[from=4-5, to=4-6]
		\arrow[from=4-6, to=3-6]
		\arrow[shorten <=3pt, Rightarrow, scaling nfold=4, no head, from=0, to=4-2]
	\end{tikzcd}\]
	\caption{Some important shapes.}
	\label{fig:shapes}
\end{figure}
}

\subsection{Extension~2: Extension Types (\rzksimpsecref{04-extension-types})}\label{ssec:extn-types}

Let $\Phi \subseteq \Psi$ be an inclusion of subshapes, in cube context $I$. An \emph{extension type} as introduced in~\cite{RiehlShulman2017}, originally due to unpublished work by Lumsdaine and Shulman, captures the strict extension of a section defined on the smaller shape $\Phi$ to the larger shape $\Psi$. Concretely, assume given a type family $I \; | \; \Psi \;|\; \Gamma \vdash A$ together with a partial section $t:I \; | \; \Phi \;|\;\Gamma \vdash a(t) : A(t) $ over the subshape $\Phi \subseteq \Psi$. Then, the corresponding \emph{extension type} has as elements the strict extensions $t:I\;|\;\Psi \;|\;\Gamma \vdash b(t) : A(t) $ such that $a|_\Phi \jdeq b$. We denote the extension type by $\exten{t:\Psi}{A(t)}{\Phi}{a}$. In case $A$ is a constant type, the ensuing extension type will be written as $\ndexten{\Psi}{A}{\Phi}{a}$.

In analogy to ordinary type-to-type function types, we can emulate shape-to-type function types by instantiating extension types by the ``empty tope'' $\varphi \defeq \bot$ and the canonical term $\tyrec_\bot$, allowing us to define the functions of shape $\Psi$ into type $A$ as $\Psi \to A \defeq \ndexten{\Psi}{A}{\bot}{\tyrec_\bot}$, and similarly for the dependent case.

\subsubsection*{Extension Extensionality (\rzkref{simplicial-hott/04-extension-types.rzk}{ExtExt})}\label{sssec:extext}

Just as in~\cite[\S 4]{RiehlShulman2017}, to make the extension types homotopically well-behaved, we also assume a version of function extensionality for extension types. In \rzk{}, we postulate an axiom that allows us to extend relative homotopies between extensions of a given partial section.

Namely, let $I$ be a cube and $\Phi \subseteq \Psi \subseteq I$ be a shape inclusion. Consider a type family $A : \Psi \to \univ$ with a partial section $a : \prod_{t:\Phi} A(t)$. As in the case of dependent functions, we may use path induction to define a map for any  $f,g : \exten{t:\Psi}{A(t)}{\Phi}{a}$ of the form
\begin{equation}\label{eq:exthtpyeq}
	 \exthtpyeq_{A,a,f,g} : (f = g) \to \exten{t:\Psi}{f(t) = g(t)}{\Phi}{\refl} \tag{\rzkref{simplicial-hott/04-extension-types.rzk}{ext-htpy-eq}}.
\end{equation}

As we did for function extensionality, we assert an extension extensionality axiom of the following form.

\begin{ax}[\rzkref{simplicial-hott/04-extension-types.rzk}{ExtExt}]\label{ax:extext}
	For any $A$, $a$, $f$, and $g$ as above, the map~\eqref{eq:exthtpyeq} is an equivalence,  \ie, there exists a term
	\[ \extext : \prod_{A, a, f, g}  \isEquiv(\exthtpyeq_{A,a,f,g}) \]
\end{ax}

In the original paper, \Cref{ax:extext} is derived from another version of the extension extensionality axiom \cite[Axiom 4.6]{RiehlShulman2017}. This version is analogous to the version of function extensionality  that states that, given a family $B:A \to \univ$, then if every fiber $Bx$ is contractible, then so is the type $\prod_{x:A} Bx$.

In the case of ordinary dependent function types, weak function extensionality is known to be equivalent to the version of function extensionality (\rzkref{hott/03-equivalences.rzk}{FunExt}). However, at the time of the writing of the paper \cite{RiehlShulman2017}, it was not known whether this equivalence also held for  extension types. Therefore, Riehl--Shulman assumed the version appearing as ~\cite[Axiom~4.6]{RiehlShulman2017}  and proved that the other desired versions, such as (\rzkref{simplicial-hott/04-extension-types.rzk}{ExtExt}), can be derived from it.

The axiom \cite[Axiom~4.6]{RiehlShulman2017}  is called \emph{relative function extensionality} (or \emph{extension extensionality}), and it reads as follows. Let $\Phi \subseteq \Psi \subseteq I$ be a shape inclusion and let $A: \Psi \to \univ$ be a family such that each $A(t)$ is contractible. Then, given $a: \prod_{t:\Phi} A(t)$, the type $\exten{t:\Psi}{A(t)}{\Phi}{a}$ is contractible. Our version (\rzkref{simplicial-hott/04-extension-types.rzk}{ExtExt}) then follows as one of the consequences established in ~\cite[Proposition~4.8]{RiehlShulman2017}.

More recently, in the process of formalizing synthetic $\infty$-category theory, Tashi Walde figured out that the standard proof that function extensionality implies weak function extensionality can be adapted to show that extension extensionality implies weak extension extensionality, and thus these axioms are logically equivalent.\footnote{\url{https://rzk-lang.github.io/sHoTT/simplicial-hott/03-extension-types.rzk/}}

\section{Synthetic \texorpdfstring{$\infty$}{∞}-Categories}\label{sec:syn-inftycat}

Simplicial type theory is a combination of the homotopical interpretation of Martin-L\"{o}f type theory with strict shapes and extension types. As demonstrated in \cite{RiehlShulman2017,BuchholtzWeinberger2023,Weinberger2022,BardomianoMartinez2022}, this framework is powerful enough to develop $\infty$-category theory synthetically, within a genuinely homotopical framework.

A category is a structure consisting of points and arrows that can be composed. To yield a synthetic notion of $\infty$-category we want to implement this idea in homotopy type theory, using the simplicial shapes.

\clearpage 
\subsection{Pre-\texorpdfstring{$\infty$}{∞}-Categories and \texorpdfstring{$\infty$}{∞}-Categories}

\subsubsection*{Hom Types (\rzkref{simplicial-hott/05-segal-types.rzk}{hom})}\label{sssec:hom-types}

\begin{figure}
\[\begin{tikzcd}
	{\Bigg\{{\small x}} & {{\small y}} & {{\small z}\Bigg\}} & {\Bigg\{{\small x}} & {{\small y}} & {{\small z}\Bigg\}}
	\arrow["\simeq"{description}, draw=none, from=1-3, to=1-4]
	\arrow["{{\small g}}", from=1-2, to=1-3]
	\arrow["{{\small f}}", from=1-4, to=1-5]
	\arrow["{{\small g}}", from=1-5, to=1-6]
	\arrow[""{name=0, anchor=center, inner sep=0}, "{{\small g\circ f}}"', curve={height=30pt}, dashed, from=1-4, to=1-6]
	\arrow["{{\small f}}", from=1-1, to=1-2]
	\arrow["{{\tiny \mathrm{comp}_{f,g}}}"{pos=0.8}, shorten <=3pt, Rightarrow, dashed, no head, from=0, to=1-5]
\end{tikzcd}\]
\caption{The Segal condition.}
	\label{fig:segal}
\end{figure}

First, we need to obtain a notion of (directed) arrow in a type $A$, which we define to be a map $f : \Delta^1 \to A$ out of the 1-simplex. The source and target of the arrow $f$ are given by the terms $f(0), f(1) : A$ respectively.

Using extension types, we can define the type of arrows with fixed source and target. For a type $A$ with elements $a,b:A$, the type of \emph{arrows} or \emph{homomorphisms} from $a$ to $b$ is
the extension type
\[ \hom_A(a,b) \defeq \ndexten{\Delta^1}{A}{\partial \Delta^1}{[a,b]}, \tag{\rzkref{simplicial-hott/05-segal-types.rzk}{hom}} \]
where $t : \partial \Delta^1 \vdash [a,b](t) : A$ is the term with $[a,b](0) \jdeq a$ and $[a,b](1) \jdeq b$.

\subsubsection*{Identity Arrows (\rzkref{simplicial-hott/05-segal-types.rzk}{id-hom})}

By the introduction rule for extension types, any element $x:A$ induces an \emph{identity arrow} $\id_x : \hom_A(x,x)$, $\id_x \defeq \lambda s.x$.

\subsubsection*{Pre-$\infty$-Categories (\rzksimpsecref{05-segal-types})}\label{sssec:segal-types}

Any type has arbitrarily complicated simplicial structure. After all, we can probe any type by an $n$-simplex for an arbitrarily high dimension $n$. When is composition of arrows, \ie, $1$-simplices defined? We want to state this as a homotopically meaningful condition: any pair of composable arrows should have a composite arrow, bounding a $2$-simplex that witnesses that the new arrow is in fact the composite of the given pair. This is classically known as the Segal condition~see \Cref{fig:cov} and~\cite{Grothendieck1961,Segal1968}. In our system, we can express it using extension types and contractibility.

\begin{defn}[pre-$\infty$-categories; \rzkrefU{simplicial-hott/05-segal-types.rzk}{is-pre-\%E2\%88\%9E-category}{is-pre-∞-category}]\label{def:segal}
	A type is a \emph{pre-$\infty$-category} or \emph{Segal type} if any composable pair of arrows has a unique composite, \ie, given a pair of arrows $f:\hom_A(x,y)$ and $g:\hom_A(y,z)$ the type of fillers \[ \sum_{h:\hom_A(x,z)} \hom_A^2(f,g;h)\]
	is contractible, where
	\[ \hom_A^2(f,g;h) \defeq \ndexten{\Delta^2}{A}{\partial \Delta^2}{[f,g;h]} \tag{\rzkref{simplicial-hott/05-segal-types.rzk}{hom2}}\] is
	the type of $2$-simplices bounded by a fixed choice of $1$-simplices:
		\begin{align*}
			 &\mathsf{is}\text{-}\mathsf{pre}\text{-}\infty\text{-}\mathsf{category}(A)  \defeq  \\ &\prod_{x,y,z:A} \prod_{f:\hom_A(x,y)} \prod_{g:\hom_A(y,z)} \!\! \isContr\left( \sum_{h:\hom_A(x,z)} \!\! \hom_A^2(f,g;h) \right)
		\end{align*}
\end{defn}

A synthetic pre-$\infty$-category is a type that admits unique composition of arrows up to contractibility.

Spelled out, this means there exists an arrow $g \circ f : \hom_A(x,z)$ acting as a composite of $g$ and $f$, together with a $2$-cell $\comp_{A,f,g} : \hom_A^2(f,g;g \circ f)$ that witnesses that the $2$-simplex bounded by $f$, $g$, and $g \circ f$, is filled, \cf~\Cref{fig:segal}. Moreover, the pair of data $g \circ f$ and $\comp_{A,f,g}$ is determined uniquely up to homotopy.

One can show that that the Segal condition of~\Cref{def:segal} can be re-expressed by saying that the type $A$ is local with respect to an inner horn inclusion.

\begin{thm}[\rzkrefU{simplicial-hott/05-segal-types.rzk}{is-pre-\%E2\%88\%9E-category-iff-is-local-horn-inclusion}{is-pre-∞-category-iff-is-local-horn-inclusion}]\label{def:segal-comp}
A type $A$ is a \emph{pre-$\infty$-category} if and only if restriction along the shape inclusion $\Lambda_1^2 \subseteq \Delta^2$ is an equivalence
\[ \isEquiv\Big(\res_A:A^{\Delta^2} \to A^{\Lambda_1^2}\Big). \tag{\rzkref{simplicial-hott/05-segal-types.rzk}{is-local-horn-inclusion}}\]
\end{thm}

\subsubsection*{Isomorphisms in Pre-$\infty$-Categories}

In a pre-$\infty$-category $A$ we may define the type of \emph{isomorphisms} by
\[ \tyiso_A(x,y) \defeq \sum_{f:\hom_A(x,y)} \isIso(f), \tag{\rzkref{simplicial-hott/10-rezk-types.rzk}{Iso}} \]
where
\[ \isIso(f) \defeq \sum_{g:\hom_A(y,x)} g \circ f = \id_x \times
\sum_{h:\hom_A(y,x)} f \circ h = \id_y.
\tag{\rzkref{simplicial-hott/10-rezk-types.rzk}{is-iso-arrow}}\]
This is in analogy with the definition of equivalences in \S\ref{sssec:hott-types}, and $\isIso(f)$ is a proposition whenever the ambient type $A$ is a pre-$\infty$-category (\rzkref{simplicial-hott/10-rezk-types.rzk}{is-prop-is-iso-arrow}).

\subsubsection*{$\infty$-Categories (\rzksimpsecref{10-rezk-types})}\label{sssec:rezk-types}

Synthetic pre-$\infty$-categories, \cf~\Cref{def:segal}, have two competing notions of sameness between terms, defined by the identity types and by the types of isomorphisms. A pre-$\infty$-category is an $\infty$-category if the notion of isomorphism just defined coincides with the notion of path in a type. This requirement captures the established notions of \emph{Rezk completeness} or \emph{local univalence}~\cite{Rezk2001,AhrensKapulkinLumsdaine2015}.

By path induction, we can define a family of comparison maps
\[ \idtoiso_A : \prod_{x,y:A} (x =_A y) \to \tyiso_A(x,y) \tag{\rzkref{simplicial-hott/10-rezk-types.rzk}{iso-eq}} \]
via path induction by
\[ \idtoiso_A(x,x,\refl_x) \defeq (\id_x,\id_x,\refl_{\id_x},\id_x,\refl_{\id_x}).\]

\begin{defn}[$\infty$-categories; \rzkrefU{simplicial-hott/10-rezk-types.rzk}{is-\%E2\%88\%9E-category}{is-∞-category}]
	A type $A$ is an \emph{$\infty$-category} or \emph{Rezk} if it is a pre-$\infty$-category and Rezk-complete:
	\begin{align*}
		\mathsf{is}&\text{-}\infty\text{-}\mathsf{category}(A)  \defeq \\ &\isSegal(A)  \times \prod_{x,y:A} \isEquiv\left(\idtoiso_{A,x,y} \right).
	\end{align*}
\end{defn}

This synthetic definition semantically translates to the well-understood notion of \emph{complete Segal} or \emph{Rezk space}, a model of $\infty$-category~\cite{Rezk2001,JoyalTierney2007,Bergner2018,Rasekh2021qc}.\footnote{More generally, this yields Rezk \emph{objects} internal to a given $\infty$-topos by~\cite{RiehlShulman2017,Shulman2019,Weinberger2022}, hence \emph{internal} $\infty$-categories~\cite{deBrito2016,Rasekh2022,Stenzel2023,MartiniWolf2023}.}

And also from an internal standpoint, the ensuing theory is quite rich. One obtains notions of functors, natural transformations, functor categories, adjunctions~\cite{RiehlShulman2017}, (co)limits~\cite{BardomianoMartinez2022}, and fibrations~\cite{RiehlShulman2017,BuchholtzWeinberger2023,Weinberger2022Sums,Weinberger2022TwoSided}, with many parallels to Riehl--Verity's \emph{$\infty$-cosmos theory}, a model-in\-de\-pend\-ent approach to $\infty$-category theory~\cite{RiehlVerity2022}.

\subsubsection*{Naturality for Free}

A useful feature of the synthetic theory is that various functoriality and naturality properties are automatically satisfied. This saves a lot of work compared to set-theoretic foundations. For instance, given pre-$\infty$-categories $A$ and $B$ any type-theoretic function $f:A \to B$ turns out to be a functor, \ie, preserves compositions and identity, up to propositional equality, see~(\rzkref{simplicial-hott/06-2cat-of-segal-types.rzk}{functors-pres-id}) and (\rzkref{simplicial-hott/06-2cat-of-segal-types.rzk}{functors-pres-comp}).\footnote{Semantically, this can be seen as the fact that any morphism between simplicial objects automatically preserves composition if they happen to be pre-$\infty$-categories.} In particular, we do not have to specify object and morphism part separately.' Similarly, given functors $f,g:A \to B$, natural transformations may be defined, using the extension types, to be arrows in the type $A \to B$, \emph{i.e.},\; $\varphi : \nat{A}{B}(f,g)$ (\rzkref{simplicial-hott/06-2cat-of-segal-types.rzk}{nat-trans}). This definition automatically yields the expected naturality squares without having to specify them, \cite[Proposition~6.6]{RiehlShulman2017}.

Further instances of automatic naturality appear in \S\ref{sssec:covariant} and \S\ref{sssec:yoneda-proof}.

\subsection{Covariant Families of $\infty$-Groupoids}

\subsubsection*{$\infty$-Groupoids (\rzksimpsecref{07-discrete})}

We are also interested in synthetic $\infty$-group\-oids, meaning $\infty$-categories where every arrow is invertible.\footnote{As shown in~\cite[\S 7]{RiehlShulman2017}, one can drop the assumption of being an $\infty$-category as it will be implied.} E.g., one can show that for any pre-$\infty$-category $A$, the hom types $\hom_A(x,y)$ are necessari\-ly $\infty$-groupoids. This matches up with the traditional theory and the intuition that $\infty$-categories are (weakly enriched) in \emph{spaces} as modeled by $\infty$-groupoids~\cite{Quillen2006}.

The groupoidal condition can be understood as a kind of \emph{discreteness condition}. To make it precise, we need a comparison of paths with arrows, similarly to our treatment of Rezk completeness, cf.~\S~\ref{sssec:rezk-types}. Namely, for a type $A$ we define
\[ \idtoarr_A : \prod_{x,y:A} (x =_A y) \to \hom_A(x,y) \tag{\rzkref{simplicial-hott/07-discrete.rzk}{hom-eq}} \]
via path induction by
\[ \idtoarr_A(x,x,\refl_x) \defeq \id_x .\]

\begin{defn}[$\infty$-groupoids; \rzkrefU{simplicial-hott/07-discrete.rzk}{is-\%E2\%88\%9E-groupoid}{is-∞-groupoid}]\label{def:disc}
	A type $A$ is an \emph{$\infty$-groupoid} or \emph{discrete} if
	\[ \isDiscrete(A) \defeq \prod_{x,y:A} \isEquiv(\idtoarr_{A,x,y}). \]
\end{defn}

This definition also yields the desired notion in the Segal object models~\cite{Bousfield1992,Stenzel2022}.

\begin{figure}
	\[\begin{tikzcd}
		{\sum_{a:A} C(a)} & u && {f_*(u)} \\
		\\
		A & x && y
		\arrow[two heads, from=1-1, to=3-1]
		\arrow["f", from=3-2, to=3-4]
		\arrow["{\mathrm{lift}_{C,f}(u)}", dashed, from=1-2, to=1-4]
	\end{tikzcd}\]
	\caption{A covariant family $C:A \to \univ$.}\label{fig:cov}
\end{figure}

\subsubsection*{Covariant Families (\rzksimpsecref{08-covariant})}\label{sssec:covariant}

The $\infty$-categorical Yo\-ne\-da lemma deals with families or \emph{fibrations} of $\infty$-groupoids indexed by a (pre-)$\infty$-category. These families $C:A \to \univ$ are supposed to be \emph{functorial} in the sense that an arrow $f : \hom_A(x,y)$ in the base $A$ should give a functor $f_* \defeq \trans_{C,f}: C(x) \to C(y)$ between the fibers.

This is achieved by the notion of \emph{covariant family}, corresponding to what semantically is often called \emph{left fibration}, after~\cite[\S 8]{Joyal2008} and \cite[\S 2.1]{Lurie2009}, see also~\cite{KazhdanVarshavsky2014,MoerdijkHeuts2015,AyalaFrancis2020,BarwickShah2018,Rasekh2023Yoneda,RiehlVerity2017,Cisinski2019}.

To define it, we have to introduce a \emph{dependent} version of the hom type, capturing arrows in the total type $\sum_{x:A}C(x)$ that get mapped to a prescribed arrow in the base. This can, once again, conveniently be formulated using extension types.

\begin{defn}[dependent hom; \rzkref{simplicial-hott/08-covariant.rzk}{dhom}]\label{def:dhom}
	Let $C:A \to \univ$ be a type family. For elements $x,y:A$, let $f:\hom_A(x,y)$ be an arrow. For elements in the fibers $u:C(x)$ and $v:C(y)$, the corresponding \emph{dependent hom type} from $u$ to $v$ is given by the extension type
	\[ \dhom_{C(f)}(u,v) \defeq \exten{t:\Delta^1}{C(f(t))}{\partial \Delta^1}{[u,v]}.\]
\end{defn}

The defining property for a covariant family $C:A \to \univ$ says that we can lift an arrow $f:\hom_A(x,y)$ in the base, given a point $u:C(x)$ in the fiber over its source, to a dependent arrow
\begin{equation*}
	\lift_{C,f,u}:\dhom_{C(f)}(u,f_*u) \tag{\rzkref{simplicial-hott/08-covariant.rzk}{covariant-transport}}
\end{equation*}
lying over $f$, and more so, uniquely up to homotopy, \cf~\Cref{fig:cov}.

\begin{defn}[Covariant family; \rzkref{simplicial-hott/08-covariant.rzk}{is-covariant}]
	Let $C:A \to \univ$ be a type family. We say $C$ is \emph{covariant} if the following proposition is inhabited:
 	\[ \prod_{x,y:A} \prod_{f:\hom_A(x,y)} \prod_{u:C(x)} \isContr\left(\sum_{v:C(y)} \dhom_{C(f)}(u,v) \right) \]
\end{defn}

As shown in~\cite[\S 8]{RiehlShulman2017}, it turns out that, over a pre-$\infty$-category $A$, covariant families $C:A \to \univ$ behave in the expected ways. Namely, the fibers are all $\infty$-groupoids~\cite[Proposition~8.18]{RiehlShulman2017}, and they are \emph{functorial} in the following sense: for elements $x,y,z:A$, morphisms $f:\hom_A(x,y)$, $g:\hom_A(y,z)$, and an element in the fiber $u:C(x)$, we get identifications
\[g_*(f_*u) = (g\circ f)_*u \; \text{and} \; (\id_x)_*u=u, \tag{\rzkref{simplicial-hott/08-covariant.rzk}{id-arr-covariant-transport}}\]
see~\cite[Proposition~8.16]{RiehlShulman2017}.

A fundamental example are the \emph{representable} covariant families of the form $\hom_A(x,-) :A \to \univ$, for $x:A$, when $A$ is a pre-$\infty$-category \\ (\rzkrefU{simplicial-hott/08-covariant.rzk}{is-covariant-representable-is-pre-\%E2\%88\%9E-category}{is-covariant-representable-is-pre-∞-category}).

Furthermore, between covariant families $C,D:A \to \univ$, a fiberwise map $\varphi : \prod_{x:A} C(x)
 \to D(x)$ is automatically \emph{natural}: for any arrow $f:\hom_A(x,y)$ and element $u:C(x)$ we have an identification
 \begin{equation}\label{eq:fam-maps-cov-is-nat}
	f_*(\varphi_x(u)) = \varphi_y(f_*u).\tag{\rzkref{simplicial-hott/08-covariant.rzk}{naturality-covariant-fiberwise-transformation}}
 \end{equation}

\section{An Overview of the \rzk~Proof Assistant}\label{sec:rzk}

Kudasov has implemented \rzk{}~\cite{Kudasov2023rzk}, the first proof assistant to support simplicial type theory.
In our work since the spring of 2023, we have been developing a library\footnote{see \url{https://emilyriehl.github.io/yoneda/}} for \Rzk,
formalizing a range of results from Riehl--Shulman's work~\cite{RiehlShulman2017},
and in addition to that also the required results from standard homotopy type theory~\cite{hottbook,Rijke2022}.
The formalizations in this paper have been written for and checked with \rzk{} \href{https://rzk-lang.github.io/rzk/en/v0.7.2/}{version 0.7.2}.

Syntax of the formalized code in \rzk{} is very close to the underlying theory,
allowing for easy correspondence between statements in the code and on paper.
However, proofs in \rzk{} may appear too detailed sometimes, since, being experimental,
\rzk{} has not yet evolved enough syntactic sugar or tools like implicit parameters,
tactics, or type classes to simplify proof construction.

In this section, we overview the key features of \rzk{} that we have relied on
in our formalization. Details about \rzk{} design and implementation are out of
the scope of this paper and should appear later in a separate paper. Still,
we should mention that the underlying implementation is a mix of general ideas
behind implementations of dependent types in Haskell~\cite{LoehMcBrideSwierstra2010}
with an experimental representation of abstract syntax with binders~\cite{Kudasov2022-arXiv-free-scoped},
and an intuitionistic sequent-based solver for the tope layer~\cite{Kudasov2023-SCAN}.

\subsection{Key Features of \rzk{}}

The kernel of \rzk{} provides the following primitive notions and capabilities.

\subsubsection*{The Universes} There are three fixed universes: \inlinerzk{CUBE} of cubes, \inlinerzk{TOPE} of topes,
    and \inlinerzk{U} of types. In \rzk{}, \inlinerzk{U}
    contains \inlinerzk{CUBE}, \inlinerzk{TOPE}, and itself, implying an unsound
    ``type in type.''\footnote{which is also present, for example, in the \libUniMath{} library~\cite{UniMath}}
    We consider such simplification acceptable for
    the time being and hope that \rzk{} will evolve proper universes in the future.

\subsubsection*{Tope Logic} This includes both cubes and topes.
    \rzk{} has built-in unit cube \inlinerzk{1} and directed interval cube \inlinerzk{2}
    (with points \inlinerzk{*₁ : 1} and \inlinerzk{0₂ : 2} and \inlinerzk{1₂ : 2} correspondingly),
    standard topes~\cite[Figure~2]{RiehlShulman2017},
    and the inequality tope \inlinerzk{s ≤ t} required for simplicial type theory.
    When done on paper, proofs in the tope logic are usually omitted as trivial,
    and we find that in our formalization project, only fairly small problems
    have been required for coherence checks. In fact, the most complicated checks we have are involved
    in the formalization of the general result for currying for extension types (\cite[\S 4.1]{RiehlShulman2017}; \rzkref{simplicial-hott/04-extension-types.rzk}{curry-uncurry}).
    \rzk{} offers full automation of the tope layer~\cite{Kudasov2023-SCAN} which helps keep
    the \rzk{} syntax and proofs simpler and automatically locate coherence issues in proof terms.

\subsubsection*{Dependent Types} \rzk{} offers  support for dependent functions  \inlinerzk{(x : A) → B x},
    dependent pairs  \inlinerzk{Σ (x : A), B x}, and identity types \inlinerzk{x =_{A} y}.
    While at the moment of writing there is no support for user-defined implicit arguments,
    identity types allow the indices to be implicit with terms \inlinerzk{x = y} and \inlinerzk{refl}
    instead of \inlinerzk{x =_{A} y} and \inlinerzk{refl_{x : A}}, resp.
    Absence of implicit arguments and full type inference in \rzk{} induces
    more explicit and verbose proof terms.

\subsubsection*{Extension Types} \rzk{} offers two separate concepts that result in support for extension types.
    First, \rzk{} allows dependent functions to have a cube or a shape (a cube restricted with a tope) argument.
    These correspond to extension types restricted to $\tyrec_\bot$ at the empty tope $\bot$.

    Second, any type is allowed to have a ``refinement,'' specifying values
    for arbitrary tope constraints.
    For example, a type \inlinerzk{A [φ ↦ x, ψ ↦ y]}
    is a refinement of type \inlinerzk{A} such that values of this type
    are \emph{computationally} equal to \inlinerzk{x} when \inlinerzk{φ} holds
    and to \inlinerzk{y} when \inlinerzk{ψ} holds.
    Of course, \inlinerzk{x} and \inlinerzk{y} must agree when \inlinerzk{(φ ∧ ψ)} holds.
    Refinements of a type are its subtypes, and \inlinerzk{A} is considered equivalent to \inlinerzk{A [⊥ ↦ recBOT]}.
    The subtyping is handled by \rzk{}, removing the need for explicit type coercions.

    Combining functions depending on shapes with such refinements yields extension types.
    For instance, $\hom_A(a,b) \defeq \ndexten{\Delta^1}{A}{\partial \Delta^1}{[a,b]}$ (\rzkref{simplicial-hott/05-segal-types.rzk}{hom}) is defined as follows:
    \begin{minted}{rzk}
#def hom (A : U) (a b : A) : U
  := (t : Δ¹) → A [t ≡ 0₂ ↦ a , t ≡ 1₂ ↦ b]
    \end{minted}

\subsubsection*{Sections and Variables} \rzk{} supports Coq-style sections,\footnote{\url{https://rzk-lang.github.io/rzk/en/v0.7.2/reference/sections.rzk/}}
allowing for locally defined assumptions (variables) which are automatically added as
parameters to definitions that use them. Importantly, \rzk{} features a mechanism
for detecting implicitly used assumptions to avoid accidental circular reasoning
in definitions. To ensure that such an implicit assumption is not accidental,
\rzk{} has the \inlinerzk{uses} syntax. For example, the Yoneda lemma (\rzkref{simplicial-hott/09-yoneda.rzk}{yoneda-lemma}) itself
is specified in a way that makes explicit the use of function
extensionality (\inlinerzk{funext}).
\begin{minted}{rzk}
#def yoneda-lemma uses (funext)
  ( A : U)
  ( is-pre-∞-category-A : is-pre-∞-category A)
  ( a : A)
  ( C : A → U)
  ( is-covariant-C : is-covariant A C)
  : is-equiv
      ((z : A) → hom A a z → C z)
      (C a)
      (evid A a C)
  := ...
\end{minted}

We find this particularly useful for readability, highlighting the use of axioms
or other assumptions (e.g. that a certain type is a pre-$\infty$-category).

\section{The \texorpdfstring{$\infty$}{∞}-Categorical Yoneda Lemma in \rzk}\label{sec:yoneda}

\subsubsection*{The Statement}

In $1$-category theory, the Yoneda lemma says the following. Given a category $\AA$ and a \emph{copresheaf\footnote{Our formalization considers the covariant case as well as the dual contravariant case.} on $\AA$}, \ie, a functor $C:\AA \to \Set$, for any $a \in \ob(\AA)$ there is a bijection
\[ \hom_{[\AA, \Set]}(\hom_{\AA}(a,-), C) \cong C(a) \]
mapping a natural transformation $\alpha$ to $\alpha(a,\id_a) \in C(a)$, naturally in both $C$ and $a$.

In the $\infty$-categorical setting, sets get replaced by $\infty$-group\-oids. Copresheaves are modeled by left fibrations \emph{aka} covariant families. Accordingly, the synthetic $\infty$-categorical Yoneda lemma reads as follows.

\begin{thm}[\rzkref{simplicial-hott/09-yoneda.rzk}{yoneda-lemma}]\label{thm:yoneda}
	Let $C:A \to \univ$ be a covariant family over a pre-$\infty$-category $A$. Then, for any $a:A$ the map
	\[ \evid_{A,a,C} : \left( \prod_{z:A} \hom_A(a,z) \to C(z) \right) \to C(a)\]
	defined by
	\[ \evid_{A,a,C}(\varphi) \defeq \varphi(a,\id_a) \tag{\rzkref{simplicial-hott/09-yoneda.rzk}{evid}} \]
	is an equivalence.
\end{thm}

Note this result holds for pre-$\infty$-categories, not just $\infty$-categories. For semantical accounts of the $\infty$-categorical Yoneda lemma see e.g.~\cite{KazhdanVarshavsky2014,Rasekh2023Yoneda,RiehlVerity2017,MartiniWolf2021}, and \cite[\S 5]{RiehlVerity2022}.

\subsubsection*{The Proof}\label{sssec:yoneda-proof}

An inverse map is constructed using the covariant transport of $C$. Namely, we define
\[ \yon_{A,a,C} :  C(a) \to  \left( \prod_{z:A} \hom_A(a,z) \to C(z) \right)  \]
by
\begin{equation}
	\label{eq:yon-natural-transformation}
	\yon_{A,a,C}(u) \defeq \lambda x.\lambda f.f_*u.
	\tag{\rzkref{simplicial-hott/09-yoneda.rzk}{yon}}
\end{equation}

In the $1$-categorical Yoneda lemma, a crucial part of the work is to show that the terms $\yon_{A,a,C}(u)$ defined by the inverse map are actually morphisms of presheaves, \ie, natural transformations. In our setting, this is, in fact, an automatic consequence from both $C$ and $\hom_A(a,-):A \to \univ$ being covariant. In the formalization considerable work goes into showing a type $A$ is a pre-$\infty$-category if and only if the type families $\hom_A(a,-)$ are covariant; for the implication relevant here, see (\rzkrefU{simplicial-hott/08-covariant.rzk}{is-covariant-representable-is-pre-\%E2\%88\%9E-category}{is-covariant-representable-is-pre-∞-category}).

	In more detail, if $A$ is a pre-$\infty$-category, $a : A$, and $C:A\to \univ$ is a covariant family, let $\varphi:\prod_{z:A} \hom_A(a,z) \to C(z)$ be a family of maps. Then for any $x,y:A$ and arrows $f:\hom_A(a,x)$ and $g:\hom_A(x,y)$, we have
\begin{equation}\label{eq:nat-cov-hom-rep}
	g_*(\varphi(x,f)) = \varphi(y, g \circ f) \footnote{(\rzkref{simplicial-hott/09-yoneda.rzk}{naturality-covariant-fiberwise-representable-transformation})}
\end{equation}
as a special case of
\[ {\eqref{eq:fam-maps-cov-is-nat}}.\]

For the Yoneda lemma, one has to show that the two composites of $ \evid_{A,a,C}$ and $\yon_{A,a,C}$ yield identities. The direction $\evid_{A,a,C} \circ \yon_{A,a,C} = \id$ is rather easy to see. Using function extensionality, we can check this pointwise. For $u:C(a)$ we have to produce an identification
\[(\lambda x.\lambda f.f_*u)(a,\id_a) = u.\]
But the left-hand side evaluates to $(\lambda x.\lambda f.f_*u)(a,\id_a) = (\id_a)_*(u)$, and the claim follows by concatenating with the path $(\id_a)_*(u) = \id_{C(a)}(u)$ given by the fact covariant transport of identities gives the identity functor of the fiber.

For the other direction, the main work is giving homotopies
\begin{align*}
    \yon_{A,a,C} ( \evid_{A,a,C}(\varphi)) (x,f)= \varphi(x,f).
\end{align*}
for all $\varphi: \prod_{z:A} \hom_A(a,z) \to C(z)$, $x:A$, and $f:\hom_A(a,x)$.

We first get a path $f_*(\varphi(a,\id_a)) = \varphi(x,f \circ \id_a)$, again by naturality \eqref{eq:nat-cov-hom-rep}. Then, using action on paths on the canonical identification $f \circ \id_a = f$ gives $\varphi(x,f \circ \id_a) = \varphi(x,f)$, and we are done.

We now have to abstract over evaluating at $x$ and $f$ using function extensionality twice, which ultimately yields the desired identification $(\yon_{A,a,C} \circ \evid_{A,a,C})(\varphi) = \varphi$, between fibered maps of type $\prod_{z:A} \hom_A(a,z) \to C(z)$. This concludes the proof of the Yoneda lemma.

\subsubsection*{The Dependent Yoneda Lemma}

The Yoneda lemma of Theorem \ref{thm:yoneda} is some sort of ``arrow-induction'' principle, but is not expressed in fully dependent form. This inspired the authors of~\cite{RiehlShulman2017} to search for a dependent generalization, proving a theorem that had not been previously known for $\infty$-categories.\footnote{There is a version of the dependent Yoneda lemma appearing in \cite[\S 5]{RiehlVerity2022}, but this was proven subsequently.}

From the \emph{dependent Yoneda lemma}, the ``absolute'' version~(\rzkref{simplicial-hott/09-yoneda.rzk}{yoneda-lemma}) can be derived.
\begin{thm}[\rzkref{simplicial-hott/09-yoneda.rzk}{dependent-yoneda-lemma}]\label{thm:dep-yoneda}
	Let $A$ be a pre-$\infty$-category, $a:A$, and
	$C: \left( \sum_{x:A} \hom_A(a,x) \right) \to \univ$ be a covariant family. Then the  map
	\[ \devid_{A,a,C} :  \left( \prod_{z:A} \prod_{f:\hom_A(a,z)} C(z,f) \right) \to C(a, \id_a)\]
	defined by
	\[ \devid_{A,a,C}(\varphi) \defeq \varphi(a,\id_a) \tag{\rzkref{simplicial-hott/09-yoneda.rzk}{dependent-evid}} \]
	is an equivalence.
\end{thm}

Note that~\Cref{thm:dep-yoneda} is reminiscent of the path induction principle, \cf~\S~\ref{sssec:id}, and hence it can be seen as a \emph{directed arrow induction principle}.

\section{Comparing the \texorpdfstring{$\infty$}{∞}- vs 1-Categorical Yoneda Lemmas}\label{sec:comparison}

A fundamental advantage of the synthetic framework for $\infty$-category is that it narrows the gap between $\infty$-category theory and 1-category theory, by moving much of the complexity inherent in homotopy coherent mathematics into the background foundation system. We can see this by comparing the formalization of the Yoneda lemma for $\infty$-categories in \Rzk\, with the formalization of the Yoneda lemma for 1-categories in other proof assistants.

Below we compare our work against 1-categorical Yoneda lemma formalizations in \libagdaunimath{} and Lean (3 and 4). Other notable formalizations of the 1-categorical Yoneda lemma appear in \libUniMath{}\footnote{\url{https://github.com/UniMath/UniMath/blob/7d7fb997dbe84b0d0107adc963281c6efb97ff60/UniMath/CategoryTheory/yoneda.v\#L325-L328}}~\cite{UniMath}, \libagdacategories{}\footnote{see \url{https://agda.github.io/agda-categories}}~\cite{HuCarette2021}, \href{https://1lab.dev}{1Lab}\footnote{\url{https://1lab.dev/Cat.Functor.Hom.html}}, Archive of Formal Proofs in Isabelle/HOL~\cite{OKeefe2005,Katovsky2010,Stark2016}, and \libmathlib{}\footnote{\url{https://leanprover-community.github.io/mathlib4_docs/Mathlib/CategoryTheory/Yoneda.html\#CategoryTheory.yoneda}}. All of these proof assistants implement some sort of dependent type theory under the hood, but the vernacular employed by the libraries \libagdacategories{}, \libmathlib{}, and the Archive of Formal Proofs is meant to reflect traditional foundations, where all types are sets, while the vernacular employed by the libraries \libagdaunimath{}, \libUniMath{}, and 1Lab is meant to reflect univalent foundations, using Voevodsky's univalence axiom to convert equivalences between types to identities. The latter perspective is inconsistent the uniqueness of identity proofs, so can only be implemented in proof assistants that support intensional identity types compatible with the univalence axiom. To capture both perspectives, the formalization discussed in \S\ref{ssec:agda-unimath} is in univalent foundations, while the formalizations discussed in \S\ref{ssec:lean-yoneda} are in traditional set-based foundations.

\subsection{1-Categorical Yoneda Lemma in \libagdaunimath{}}\label{ssec:agda-unimath}
As part of this project, we contributed a formalization of the Yoneda lemma for precategories to the \libagdaunimath{} library,\footnote{\url{https://unimath.github.io/agda-unimath/category-theory.yoneda-lemma-precategories.html}} which describes itself as a ``community-driven effort aimed at formalizing mathematics from a univalent point of view.'' This library contains notions of \emph{precategories} and \emph{categories}, which parallel our pre-$\infty$-categories and $\infty$-categories, except their hom-types are sets, as is appropriate to 1-category theory. Both proofs follow the same outline, proving that (\rzkref{simplicial-hott/09-yoneda.rzk}{evid}) is an equivalence by constructing a two-sided inverse. A point of difference in the \libagdaunimath{} proof is that the data of the inverse involves both the function \eqref{eq:yon-natural-transformation} together with a proof of its naturality. As with our proof in \Rzk{}, one of the composites is directly identifiable with the identity, while the other requires a calculation together with two instances of function extensionality.

Other differences arise from the varying ways that categorical data is encoded in \Rzk{} vs \libagdaunimath{}. There, precategories are types with additional structure while here pre-$\infty$-categories are types satisfying a property. There, representables are encoded as functors valued in the precategory of sets, while here representables are encoded as covariant type families. These differences have more of an effect on the syntax of the proof than its structural content.

\subsection{1-Categorical Yoneda Lemma in Lean}\label{ssec:lean-yoneda}
At our request, Sina Hazratpour wrote a Lean formalization of the 1-categorical Yoneda lemma, first as a self-contained formalization in Lean~3,\footnote{\url{https://github.com/sinhp/CovariantYonedaLean3}} with the proof of the Yoneda lemma later updated to Lean~4.\footnote{\url{https://github.com/sinhp/CovariantYonedaLean4}} Formal proofs in Lean are quite different than formal proofs in \Rzk{} or in Agda because of the use of automation tactics in the interactive theorem proving mode, allowing the user to rewrite along known identifications or ``simplify'' the goal using known lemmas. In addition, Lean's use of type classes and automatic instance inference simplifies the syntax in the statement of the Yoneda lemma, as compared with the \libagdaunimath{} proof.

In the Lean~3 proof, the naturality of \eqref{eq:yon-natural-transformation} must again be checked explicitly via a proof that involves unfolding the definition of the representable functor and using the fact that functors preserve composition. The remainder of the proof proceeds as before. Interestingly, in the Lean~4 proof, Hazratpour proves a lemma --- \eqref{eq:nat-cov-hom-rep} in the case where $f$ is $\id_a$ --- and then feeds it to the tactic \texttt{aesop\_cat},\footnote{Aesop (Automated Extensible Search for Obvious Proofs) is a proof search tactic for Lean~4; see \url{https://github.com/JLimperg/aesop}} which then automatically verifies the naturality of \eqref{eq:yon-natural-transformation} and checks that the Yoneda maps are inverses.

\section{Conclusions and Future Work}\label{sec:future}

We hope that the \Rzk{} proof assistant will provide a tool that may make $\infty$-category theory easier to learn. To that end we invite new collaborators to help us formalize other results from $\infty$-category theory. Indeed, some of this work is already underway in a new repository,\footnote{\url{https://github.com/rzk-lang/sHoTT}} which originated as a clone of our Yoneda repository, as we now describe.

\subsection{Adjunctions}

One application of the dependent Yoneda lemma is to the theory of adjunctions between $\infty$-categories. The standard logically equivalent definitions of an adjunction are encoded by various types that define a \emph{transposing adjunction}, \emph{half-adjoint diagrammatic adjunction}, or a \emph{bidiagrammatic adjunction}. These definitions have been formalized\footnote{\url{https://rzk-lang.github.io/sHoTT/simplicial-hott/11-adjunctions.rzk/}} and we have begun the lengthy task of formally proving the equivalences established in \cite[\S 11]{RiehlShulman2017}.

\subsection{Limits and Colimits}

In \cite{BardomianoMartinez2022}, Bardomiano Mart\'{i}nez introduces limits and colimits of diagrams valued in pre-$\infty$-categories and proves that right adjoints between Segal types preserve limits. Bardomiano Mart\'{i}nez has formalized these definitions\footnote{\url{https://rzk-lang.github.io/sHoTT/simplicial-hott/13-limits.rzk/}} and plans to work with us to formalize his results. Once this is done, we would like to explore further developments of the theory of limits and colimits.

\subsection{The Cocartesian $\infty$-Categorical Yoneda Lemma}

As a future endeavor, it is desirable to extend the existing fibrational theory from $\infty$-groupoid-valued functorial families to $\infty$-category-valued functorial families. Building on the semantical theory of (co)cartesian fibrations~\cite{Joyal2008,Lurie2009,AyalaFrancis2020,BarwickShah2018,RiehlVerity2017,RiehlVerity2022,Rasekh2023Fib,Martini2022} these so-called \emph{(co)cartesian families} have been studied in the synthetic setting of simplicial type theory in~\cite{BuchholtzWeinberger2023,Weinberger2022PhD}. Cocartesian fibrations play a crucial role when computing limits of $\infty$-categories~\cite[Subsection~3.3.3]{Lurie2009}, studying moduli spaces in geometry~\cite[\S 1]{Lurie2018}, and higher algebraic objects such as symmetric monoidal $\infty$-categories~\cite[2.1.2.13]{Lurie2017}. A central theorem is the (co-)cartesian Yoneda lemma, \cf~\cite[\S 5.7]{RiehlVerity2022}. It reads similar to the version of the discrete Yoneda Lemma discussed in~\S\ref{sec:yoneda}: given a \emph{cocartesian} family $C : A \to \univ$ over an $\infty$-category $A$ (all of whose fibers are also $\infty$-categories), we want to classify families of functions $\prod_{x:A} \hom_A(a,x) \to C(x)$, for a fixed element $a:A$. However, we will have to restrict to the type $\prod_{\pair{x}{f}:\sum_{y:A}  \hom_A(a,y)}^{\mathrm{cocart}} C(x)$ of sections $\sigma$ such that for every $x':A$, $f':\hom_A(a,x')$, and $g : \hom_A(x,x')$ with $g \circ f = f'$ the dependent morphism $\sigma(g)$ is what is called a \emph{cocartesian arrow}, \ie, satisfies a certain initiality property (\rzkref{simplicial-hott/12-cocartesian.rzk}{is-cocartesian-arrow}). A family $C : A \to \univ$ whose total type is an $\infty$-category (as well as all its fibers) is cocartesian exactly if it admits enough lifts of arrows in the base $A$ to a cocartesian dependent arrow (\rzkref{simplicial-hott/12-cocartesian.rzk}{is-cocartesian-family}).

The cocartesian Yoneda lemma then states that the evaluation map
\[ \evid^C_a : \left( \prod_{\pair{x}{f}:\sum_{y:A}  \hom_A(a,y)}^{\mathrm{cocart}} C(x)  \right) \to C(a) \]
is an equivalence. The \emph{dependent version} of this theorem, semantically established by Riehl--Verity in~\cite[5.7.2]{RiehlVerity2022}, can, again, be seen as a (one-sided) \emph{directed arrow induction principle}, analogous to the well-known path induction principle for the identity types in standard Martin-L\"{o}f~type~theory.

Efforts in the direction of formalizing Buchholtz--Wein\-berger's proof of the cocartesian Yoneda lemma from~\cite[\S 7]{BuchholtzWeinberger2023} in~\rzk\,are under way, but will require formalizing if not all then at least some of the preliminary structural properties and operations for cocartesian families from~\cite[\S 5]{BuchholtzWeinberger2023}.

\subsection{Improvements to \Rzk}

We note a few improvements for \rzk{} that would positively affect this
and future formalization projects. First, supporting term inference
and implicit arguments would help reduce the size of formalizations and,
consequently, assist with readability. Second, the current implementation
lacks incremental typechecking and proper module support, which makes the feedback
on changes less immediate. Finally, while a minimal integration with an IDE
exists,\footnote{there is a VS Code extension for \rzk{} at \url{https://github.com/rzk-lang/vscode-rzk}}
it still has to acquire proper language server support.
We note also that \rzk{}'s experimental
diagram rendering feature\footnote{\url{https://rzk-lang.github.io/rzk/en/v0.7.2/reference/render.rzk/}} (which is useful on small examples)
could be extended further to assist with visualizations
(or even interactive capabilities) for statements and constructions in simplicial type theory.

\subsection{Extensions of Simplicial Type Theory}

The simplicial type theory is not sufficiently powerful to prove all results of $\infty$-category theory contained for instance in \cite{Lurie2009}. A longer range goal would be to further extend this synthetic framework by including directed higher inductive types to freely generate $\infty$-categories, universes to classify covariant fibrations and cocartesian fibrations, and modalities for opposite $\infty$-categories and the $\infty$-groupoid core as outlined in \cite{BuchholtzWeinberger2019}; see also~\cite{Shulman2018,MyersRiley2023,North2019,GratzerKavvosNuytsBirkedal2020,AagaardKristensenGratzerBirkedal2022,LicataShulmanRiley2017,LicataOrtonPittsSpitters2018}. If such theoretical developments were paired with experimental extensions to \Rzk, that would greatly aid the process of exploring the expanded formal system.

\begin{acks}
We are very grateful to Benedikt Ahrens, who first suggested the project of creating a proof assistant for the simplicial type theory. Fredrik Bakke contributed formalizations concerning the 2-category of synthetic pre-$\infty$-categories and made invaluable improvements to the professionalization of the repository, drafting a style guide, overseeing its implementation, and suggesting improvements to our github workflow. Sina Hazratpour produced a formalized proof of the 1-categorical Yoneda lemma in Lean to provide a useful direct comparison. Abdelrahman Abounegm has contributed a \rzk{} plugin\footnote{\url{https://github.com/rzk-lang/mkdocs-plugin-rzk}} for \href{https://www.mkdocs.org}{\texttt{MkDocs}} allowing for hyperlinks to the syntax-highlighted code used in this paper. He has also worked on the VS Code extension for Rzk with an experimental auto-formatting feature, improving our experience and helping us enforce the style guide.
The comments from the anonymous reviewers greatly improved the published paper.
Finally, we are grateful to Denis-Charles Cisinski, Clara L\"{o}h, and Philipp R\"{u}mmer, the organizers of the ``Interactions of Proof Assistants and Mathematics'' Summer School in Regensburg, for giving us a venue to present our work and recruit new collaborators to the project of formalizing synthetic $\infty$-category theory.

\end{acks}

\bibliographystyle{ACM-Reference-Format}
\bibliography{yoneda-rzk}

\end{document}